\documentclass[10pt,a4paper,twoside]{article}
\usepackage[T1]{fontenc}
\usepackage{amsmath,amssymb,amsthm,mathtools}
\usepackage{lmodern}
\usepackage{fancyhdr}
\usepackage{microtype}
\usepackage[left=40mm,right=40mm,top=36mm,bottom=43mm,headheight=13pt,headsep=5mm]{geometry}
\usepackage[hidelinks,pdfauthor={SHOOT KOEBISU},pdftitle={Principal angles between annihilators of real sedenion zero divisors}]{hyperref}

\makeatletter
\renewcommand{\@maketitle}{%
  {\centering\LARGE\mdseries\@title\par}
  \ifx\@author\@empty\else\vspace{8pt}{\centering\small\begin{tabular}[t]{@{}c@{}}\@author\end{tabular}\par}\fi
  \vspace{24pt}}
\makeatother
\renewenvironment{abstract}{\begin{quotation}\small\noindent\textbf{Abstract}\par\smallskip\noindent\ignorespaces}
  {\end{quotation}\vspace{5pt}}
\newtheoremstyle{displayresult}{\topsep}{\topsep}{\itshape}{}{\bfseries}{.}{\newline}{}
\theoremstyle{displayresult}
\newtheorem{theorem}{Theorem}
\theoremstyle{plain}
\newtheorem{lemma}[theorem]{Lemma}
\newtheorem{proposition}[theorem]{Proposition}
\newtheorem{corollary}[theorem]{Corollary}
\newcommand{\R}{\mathbb R}
\newcommand{\C}{\mathbb C}
\newcommand{\Oct}{\mathbb O}
\newcommand{\Sed}{\mathbb S}
\newcommand{\im}{\operatorname{Im}}
\newcommand{\spn}{\operatorname{span}}
\newcommand{\ip}[2]{\langle #1,#2\rangle}
\DeclareMathOperator{\tr}{tr}
\DeclareMathOperator{\diag}{diag}
\DeclareMathOperator{\dist}{dist}
\title{Principal angles between annihilators\\of real sedenion zero divisors}
\author{SHOOT KOEBISU}
\date{}
\begin{document}
\maketitle
\thispagestyle{empty}
\begin{abstract}
Every nonzero zero divisor of the real sedenions annihilates a real four-plane. We determine the possible principal angles between two such planes with prescribed octonionic product axes. If the axes make angle $\theta$ and $k=\cos(\theta/2)$, the four principal cosines are $(p,p,q,q)$, where $k^2\le p\le k$ and $0\le q\le p$. Every point of this region occurs for each fixed pair of axes. For a specified pair of zero divisors we give a quadratic formula in two mixed invariants. An explicit deformation keeps both quaternion subalgebras and the angle between the axes fixed while changing the principal angles of the annihilators. Metric projection onto the zero-divisor cone transfers the formulas to least singular subspaces of left multiplication.
\end{abstract}

\section{The fixed-axis problem}
Write a nonzero zero divisor as $z=(u,v)$ in $\Sed=\Oct\oplus\Oct$.
Its components are orthogonal imaginary octonions of equal length. Their product determines a unit imaginary octonion $c_z=uv/|uv|$, which we call the \emph{axis} of $z$. Let $K_z=\{y\in\Sed:zy=0\}$. For a unit zero divisor, left multiplication $L_z$ has squared singular values $0,1,2$ of multiplicities $4,8,4$, respectively \cite[Prop.~7.4]{biss}. Thus the singular values of each factor separately give no information about the relative position of $K_z$ and $K_w$.

For two real four-planes we order their principal angles increasingly in $[0,\pi/2]$. Their cosines are the singular values of the matrix of inner products of two orthonormal bases \cite{bg}. They record, in particular, the intersection dimension (the multiplicity of $1$) and orthogonality (all cosines zero).

\begin{theorem}[Exact range at fixed axes]\label{th:range}
Fix unit $c,d\in\im\Oct$, and set $k=\sqrt{(1+\ip cd)/2}$.
As $z,w$ vary among unit zero divisors with $c_z=c$ and $c_w=d$, the ordered principal cosines of $K_z,K_w$ range over precisely
\begin{equation}
 \big\{(p,p,q,q):\ k^2\le p\le k,\quad 0\le q\le p\big\}.
 \label{eq:region}
\end{equation}
\end{theorem}
In particular, orthogonality of the annihilators is equivalent to the axes being antipodal. At perpendicular axes, the smallest principal angle ranges from $\pi/4$ to $\pi/3$, although the largest can reach $\pi/2$.

\section{Annihilators and their complex structures}

We use the doubling convention
\[
 (a,b)(x,y)=(ax-\bar y b,\,ya+b\bar x)
\]
at each Cayley--Dickson step, and the Euclidean inner product on both summands. Biss, Christensen, Dugger and Isaksen identify the eigenspaces of $L_z^*L_z$ \cite[Thm.~7.2, Prop.~7.4]{biss}. Reggiani studies the Riemannian geometry of the normalized zero-divisor locus and of normalized pairs satisfying $zw=0$ \cite{reggiani}. Guterman and Zhilina describe the connected components of the zero-product graph using lines in $\im\Oct$ \cite{gz}. Our pair $z,w$ is arbitrary, and the oriented axes distinguish equal from antipodal axes. Principal angles of associative planes in $\im\Oct$ enter the work of Jost, Xin and Yang \cite[Sec.~3]{jxy}. Section~\ref{sec:mixed} gives a family in which those planes themselves stay fixed.

Let $H_z$ be the quaternion subalgebra generated by $u,v$. For a unit $c\in\im\Oct$, write $J_ch=ch$. This orthogonal complex structure preserves $H_z$ and $H_z^\perp$ when $c=c_z$. Define the isometry
\[
 F_c(h)=2^{-1/2}(h,-J_ch),\qquad h\in\Oct.
\]
The graph formula of \cite[Prop.~7.4]{biss}, unchanged by real scaling of $z$, reads
\begin{equation}
 K_z=F_{c_z}(H_z^\perp).
 \label{eq:graph}
\end{equation}
On $(\Oct,J_c)$ we use the Hermitian inner product
$[h,g]_c=\ip hg+\mathrm i\ip{J_ch}g$, linear in $g$.
Set $V_c=\spn_{\R}\{1,c\}^\perp$. This is the complex orthogonal complement of the complex line generated by $1$.

\begin{lemma}\label{lem:flags}
For fixed $c$, the spaces $H_z^\perp$ with $c_z=c$ range over all complex two-planes in the complex three-space $V_c$.
\end{lemma}
\begin{proof}
A complex two-plane in $V_c$ has complex orthogonal complement there generated by some unit $a\in V_c$. The elements $c,a$ generate the associative subalgebra $H=\spn_\R\{1,c,a,ca\}$. The element $z=(a,ca)/\sqrt2$ is a unit zero divisor, and $a(ca)=c$ in $H$. Thus $H_z=H$ has the prescribed orthogonal complement. The converse follows from $1,c\in H_z$ and $J_cH_z=H_z$.
\end{proof}

The comparison map between the two ambient graphs is
\begin{equation}
 T=F_c^*F_d=\tfrac12(I-J_cJ_d),\qquad
 TJ_d=J_cT,\qquad T^*T=k^2I.
 \label{eq:T}
\end{equation}
Indeed $J_c^*=-J_c$, and polarization of $J_x^2=-|x|^2 I$ for imaginary $x$ gives $J_cJ_d+J_dJ_c=-2\ip cd I$. In particular, $J_cJ_d$ acts by $h\mapsto c(dh)$, with the indicated bracketing. For $k>0$, $U=T/k$ is a unitary map from $(\Oct,J_d)$ to $(\Oct,J_c)$.

\section{Proof of the range theorem}

Let $P_M$ denote real orthogonal projection onto $M$. It is also Hermitian projection whenever $M$ is a complex subspace. Write $A=H_z^\perp$, $B=H_w^\perp$. The real principal cosines sought are the singular values of $P_AT|_B$, by \eqref{eq:graph}. This map is complex linear, so each complex singular value occurs twice over $\R$. If $k=0$, \eqref{eq:T} gives the singleton in \eqref{eq:region}. Assume henceforth that $k>0$.

\begin{lemma}\label{lem:three}
The complex singular values of $D=P_{V_c}U|_{V_d}$ are $1,1,k$.
\end{lemma}
\begin{proof}
In the complex four-space $(\Oct,J_c)$ the hyperplanes $V_c$ and $UV_d$ have unit normals $1$ and $U1$. Since $T1=(1-cd)/2$, its real part is $k^2$ and its component along $c$ is zero. Thus $[1,U1]_c=k$. Projection between two complex hyperplanes is isometric on their common part. On the complex two-space spanned by independent unit normals, the absolute inner product of their perpendicular unit directions equals that of the normals. This gives $1,1,k$ for $k<1$. For $k=1$ the hyperplanes coincide.
\end{proof}

\begin{proof}[Proof of Theorem~\ref{th:range}]
By Lemma~\ref{lem:flags} the two complex planes $A\subset V_c$ and $B\subset V_d$ can be chosen independently. By Lemma~\ref{lem:three} we may use singular bases and identify $D$ with $\diag(1,1,k)$ on $\C^3$. The required complex cosines are $k$ times the singular values of $P_AD|_B$.

Both latter singular values are at most $1$. Since $D$ is invertible, the complex two-planes $B$ and $D^{-1}A$ meet in a line or a plane. For any unit $x$ in their intersection, $\|P_ADx\|=\|Dx\|\ge k$. Thus the larger singular value is at least $k$.

Conversely prescribe $k\le s\le1$ and $0\le r\le s$. Take a unit $v\in\spn_\C\{e_2,e_3\}$ with $\|Dv\|=s$, possible by varying continuously between $e_2$ and $e_3$. Take $B=\spn_\C\{e_1,v\}$. In the target put $f=Dv/s$ and choose a unit $g$ perpendicular to both $e_1$ and $f$. Set
\[
 A=\spn_\C\{f,\ r e_1+\sqrt{1-r^2}\,g\}.
\]
In the indicated orthonormal bases, the projected images of $v,e_1$ are orthogonal with lengths $s,r$. Multiplication by $k$ gives $(p,q)=(ks,kr)$ throughout \eqref{eq:region}, including its boundary.
\end{proof}

This proof is a two-plane instance of singular-value interlacing for submatrices \cite{thompson}. The octonionic input identifies the three ambient singular values and realizes every source and target plane by actual unit sedenion zero divisors.

\section{Recovering the angles of a specified pair}\label{sec:formula}
The following formula uses vectors from the input quaternion subalgebras. It requires only two complex inner products per column. Choose unit $a\in H_z\cap V_c$, $b\in H_w\cap V_d$, with $c=c_z$, $d=c_w$. For example, one may take the normalized first octonionic component of each zero divisor. Put $(e_1,e_2)=(1,a)$ and $(f_1,f_2)=(1,b)$. Define
\begin{equation}
 M_{ij}=[e_i,Tf_j]_c,\qquad
 S=\sum_{i,j=1}^2 |M_{ij}|^2,\qquad P=|\det M|^2 .
 \label{eq:invariants}
\end{equation}

The entries can be evaluated directly from marked quaternion frames. For $e=e_i$ and $f=f_j$, expansion of \eqref{eq:T} gives
\[
 2M_{ij}=\ip ef+\ip{ce}{df}
       +\mathrm i\big(\ip{ce}f-\ip e{df}\big).
\]
In particular $M_{11}=k^2$. The other entries record inner products involving the two marked frames.

\begin{theorem}[Pair formula]\label{th:pair}
For a specified pair the two principal cosines, counted once per complex dimension, are given by
\begin{equation}
 p^2=\frac{S+\sqrt{S^2-4P}}2,\qquad
 q^2=\frac{S-\sqrt{S^2-4P}}2.
 \label{eq:quadratic}
\end{equation}
The scalars $S,P$ are independent of the choices of $a,b$ and invariant under simultaneous octonion automorphisms.
\end{theorem}
\begin{proof}
The displayed bases are unitary for $H_z,H_w$, so $M$ represents $P_{H_z}T|_{H_w}$. If $k>0$, decompose the unitary $U$ into square blocks with respect to $H_w\oplus H_w^\perp$ and $H_z\oplus H_z^\perp$, and write these blocks as $U_{ij}$. Unitarity gives
\[
 U_{11}U_{11}^*=I-U_{12}U_{12}^*,\qquad
 U_{22}^*U_{22}=I-U_{12}^*U_{12}.
\]
The right sides have the same spectrum since $U_{12}$ is square. Hence the singular values of $M=kU_{11}$ agree with those of $kU_{22}$, which determine the annihilator angles. For $k=0$ both blocks of $T$ vanish. The characteristic polynomial of $M^*M$ is $X^2-SX+P$. Unitary basis changes preserve both. Octonion automorphisms preserve multiplication and the real inner product and hence preserve every construction in \eqref{eq:invariants}.
\end{proof}

\begin{corollary}[Equal subalgebras and equal axes]\label{cor:special}
If $H_z=H_w$, all four principal angles equal $\tfrac12\arccos\ip cd$.
If $c=d$, the cosines are $(1,1,m,m)$ with $m=|[a,b]_c|$.
In particular the intersection dimension equals $4$ for identical axis and subalgebra, $2$ for identical axis and distinct subalgebras, and $0$ for distinct axes.
\end{corollary}
\begin{proof}
In the first case the common orthogonal complement is preserved by both complex structures and hence by $T$. Equation~\eqref{eq:T} gives four cosines $k$. In the second case $T=I$ and $M=\diag(1,[a,b]_c)$. Here $m=1$ exactly when $a,b$ generate the same complex line in $V_c$. For distinct axes $p\le k<1$.
\end{proof}

\section{A deformation with fixed quaternion subalgebras}\label{sec:mixed}

How much of \eqref{eq:invariants} can be read from the unmarked subalgebras and the axial angle? Even their joint knowledge leaves a continuous ambiguity.

\begin{proposition}\label{prop:family}
There exist fixed $H,H'\subset\Oct$ and axes $c_t\in H$, $d\in H'$ with $c_t\perp d$ for $0\le t\le\pi/2$, for which
\begin{equation}
 p^2=(1+\cos t)/4,\qquad q^2=(1-\cos t)/4 .
 \label{eq:family}
\end{equation}
In this family $S=1/2$ throughout and $P=\sin^2(t)/16$.
\end{proposition}
\begin{proof}
Double the quaternions with $ij=k_0$ to form $\Oct=\mathbb H\oplus\mathbb H\ell$, using the multiplication of Section~2 and $h\ell=(0,h)$. Choose
\[
 H=\spn_\R\{1,i,j,k_0\},\quad H'=\spn_\R\{1,i,\ell,i\ell\},\quad
 c_t=i\cos t+j\sin t,\quad d=\ell .
\]
Lemma~\ref{lem:flags} realizes these marked subalgebras by unit zero divisors. In the ordered real bases $(\ell,i\ell,j\ell,k_0\ell)$ of $H^\perp$ and $(j,k_0,j\ell,k_0\ell)$ of $(H')^\perp$, doubling gives the matrix
\[
 2P_{H^\perp}T|_{(H')^\perp}=
 \begin{pmatrix}
 -\sin t&0&0&0\\
 0&-\sin t&0&0\\
 0&\cos t&1&0\\
 -\cos t&0&0&1
 \end{pmatrix}.
\]
For instance $J_dj=-j\ell$, so $2Tj=j-(\sin t)\ell-(\cos t)k_0\ell$. After projection this gives the first column. Grouping columns $(1,4)$ and $(2,3)$ gives two orthogonal blocks. After dividing the displayed matrix by two, each block has squared singular values with sum $1/2$ and product $\sin^2(t)/16$, giving \eqref{eq:family}.
\end{proof}

At $t=0$ the real angles are $(\pi/4,\pi/4,\pi/2,\pi/2)$. At $t=\pi/2$ they are all $\pi/3$. The subalgebras meet in $\spn_\R\{1,i\}$ throughout. The moving axis changes its position relative to this common imaginary line. In particular the axial angle and all principal angles between the quaternion subalgebras together fail to determine the annihilator angles.

\section{Least singular directions near the zero-divisor cone}

Let $Z$ be the cone of zero divisors with its vertex included. For $x=(a_0+u,b_0+v)$, where $a_0,b_0\in\R$ and $u,v\in\im\Oct$, let $B_x=[u\ v]$ in an orthonormal basis of $\im\Oct$. Write $E_x$ for the eigenspace of $L_x^*L_x$ at its smallest eigenvalue.

\begin{proposition}\label{prop:nearest}
If $B_x$ has rank two, $x$ has a unique Euclidean nearest point $z(x)$ in $Z$. Its real components vanish and its imaginary-component matrix is
\[
 B_{z(x)}=\tfrac12(s_1+s_2) B_x(B_x^{\mathsf T}B_x)^{-1/2},
\]
where $s_1\ge s_2>0$ are the singular values of $B_x$. Moreover
\[
 E_x=K_{z(x)},\qquad
 \dist(x,Z)^2=a_0^2+b_0^2+\tfrac12(s_1-s_2)^2 .
\]
\end{proposition}
\begin{proof}
The imaginary matrices in $Z$ have the form $tQ$, with $Q^{\mathsf T}Q=I_2$ and $t\ge0$. In minimizing $\|B_x-tQ\|_F^2$, the maximal value of $\tr(Q^{\mathsf T}B_x)$ is $s_1+s_2$. Indeed, for a thin singular value decomposition $B_x=U\diag(s_1,s_2)V^{\mathsf T}$ it is a positively weighted sum of the inner products of corresponding unit columns of $QV$ and $U$. Equality requires $QV=U$, which identifies the unique optimal frame. Minimizing the remaining quadratic gives $t=(s_1+s_2)/2$ and the stated distance.

By \cite[Prop.~3.10 and Thm.~7.2]{biss}, $E_x=F_c(H^\perp)$, where $H$ is generated by $u,v$ and $c=\im(uv)/|\im(uv)|$. The sign agrees with \eqref{eq:graph} since the cross product in \cite{biss} is $\im(u\bar v)$. Right multiplication of $B_x$ by a positive definite two-by-two matrix preserves its oriented column plane. It therefore preserves $H$ and $c$, since $(u,v)\mapsto\im(uv)$ is alternating bilinear on imaginary inputs. This gives the same graph for $z(x)$.
\end{proof}

Consequently, for two rank-two inputs $x,y$, Theorem~\ref{th:pair} computes the principal angles of $E_x,E_y$ from their imaginary columns. As an example, take any zero divisors from Proposition~\ref{prop:family} and add a nonzero real scalar to the first octonionic component of each. Left multiplication then has strictly positive least singular value by \cite[Prop.~3.10]{biss}, while the least singular directions still follow the full variation in \eqref{eq:family}. For inputs of rank at most one, all singular values coincide \cite[Props.~3.10 and 7.1]{biss}. The least singular space then has dimension sixteen.

\end{document}